\documentclass[aap,MSNbibl,citesort,dvips]{arximspdf}
\usepackage{mathrsfs}
\usepackage{mathbh}
\usepackage{accents}
\usepackage{graphicx}

\doi{10.1214/12-AAP852} 
\volume{23}
\issue{2}
\pubyear{2013}
\firstpage{693}
\lastpage{721}

\makeatletter

\renewcommand{\underbar}{\underaccent{\bar}}

\newcommand{\eqref}[1]{(\ref{#1})}

\newtheorem{theorem}{Theorem}

\newproclaim{definition}{Definition}[section]

\newtheorem{lemma}[definition]{Lemma}

\newproclaim{remark}[definition]{Remark}

\newcommand{\cal}{\mathcal}

\newcommand{\ecal}{\mathcal{E}}
\newcommand{\scal}{\mathcal{S}}
\newcommand{\diam}{\operatorname{Diam}}

\newcommand{\T}{{\mathbf{T}}}

\newcommand{\eps}{\varepsilon}

\newcommand{\ind}{\mathbh{1}}

\newcommand{\prob}{\mathbb{P}}
\newcommand{\expec}{\mathbb{E}}
\newcommand{\var}{\operatorname{Var}}

\newcommand{\poly}{\operatorname{poly}}

\newcommand{\birth}{\lambda}
\newcommand{\death}{\mu}
\newcommand{\mutate}{\eta}
\newcommand{\s}{\sigma}
\newcommand{\dist}{\mathcal{D}}
\newcommand{\edist}{\widehat{\dist}}
\newcommand{\epartdist}{\widetilde{\dist}}
\renewcommand{\path}{\mathrm{P}}
\newcommand{\hamming}{\widehat{\mathcal{H}}}
\newcommand{\deviate}{\Delta}
\newcommand{\maxt}{\bar t}

\newcommand{\A}{\mathtt{A}}
\newcommand{\G}{\mathtt{G}}
\newcommand{\C}{\mathtt{C}}
\renewcommand{\T}{\mathtt{T}}
\newcommand{\gap}{-}
\newcommand{\data}{\mathbf{D}}
\newcommand{\good}{\ecal}

\newcommand{\weight}{\omega}
\newcommand{\netrate}{\delta}
\newcommand{\totalrate}{\phi}

\newcommand{\js}[1]{S_{#1}}
\newcommand{\realiz}{\mathcal{I}}
\newcommand{\surv}{\mathrm{S}}
\newcommand{\nsurv}{{\mathrm{NS}}}

\makeatother

\begin{document}
\begin{frontmatter}

\title{Alignment-free phylogenetic reconstruction:
Sample complexity via a branching process analysis}
\runtitle{Alignment-free phylogenetic reconstruction}

\begin{aug}
\author[A]{\fnms{Constantinos} \snm{Daskalakis}\thanksref{t1}\ead[label=e1]{costis@csail.mit.edu}}
\and
\author[B]{\fnms{Sebastien} \snm{Roch}\corref{}\thanksref{t2}\ead[label=e2]{roch@math.ucla.edu}}
\runauthor{C. Daskalakis and S. Roch}
\affiliation{Massachusetts Institute of Technology and\break University of
California, Los~Angeles}
\address[A]{Department of Electrical Engineering\\
\quad and Computer Science\\
Massachusetts Institute of Technology\\
32 Vassar Street\\
Cambridge, Massachusetts 02139\\
USA\\
\printead{e1}}
\address[B]{Department of Mathematics\\
\quad and Bioinformatics Program\\
University of California, Los Angeles\\
520 Portola Plaza\\
Los Angeles, California 90095-1555\\
USA\\
\printead{e2}} 
\end{aug}

\thankstext{t1}{Supported by a Sloan Foundation Fellowship and NSF Award
CCF-0953960 (CAREER) and CCF-1101491.}

\thankstext{t2}{Supported by NSF Grant DMS-10-07144.}

\received{\smonth{9} \syear{2011}}
\revised{\smonth{2} \syear{2012}}

%
\begin{abstract}
We present an efficient phylogenetic reconstruction algorithm allowing
insertions and deletions which provably achieves a sequence-length
requirement (or sample complexity) growing polynomially in the number
of taxa. Our algorithm is distance-based, that is, it relies on
pairwise sequence comparisons. More importantly, our approach largely
bypasses the difficult problem of multiple sequence alignment.
\end{abstract}

%

\end{frontmatter}

\section{Introduction}

We introduce a new efficient algorithm for the
\textit{phylogenetic tree reconstruction} (PTR) problem
which rigorously accounts for insertions and deletions.\vspace*{10pt}

\textit{Phylogenetic background.}\quad
A \textit{phylogenetic tree} or \textit{phylogeny}
is a tree representing the speciation history of a group of organisms.
The leaves of the tree are
typically existing species.
The root corresponds to their \textit{most recent
common ancestor} (MRCA).
Each branching in the tree indicates a speciation event.
It is common to assume that
DNA evolves according to a Markovian substitution process on this phylogeny.
Under such a model,
a \textit{gene} is a sequence in $\{\A,\G,\C,\T\}^k$.
Along each edge of the tree, each site independently mutates
according to a Markov rate matrix.
The length of a branch
is a measure of the amount of substitution along that branch.
The precise definition of a branch length depends on the model of evolution.
For roughly constant mutation rates, one can think of the branch length
as proportional to the amount of time elapsed along a branch.
The PTR problem
consists of estimating a phylogeny from the genes observed at
its leaves.
We denote the leaves of a tree by $[n] = \{1,\ldots,n\}$
and their sequences by $\s_1,\ldots,\s_n$.

The model of sequence evolution above is simplistic:
it ignores many mutational events that DNA undergoes
through evolution. At the gene level, the most important
omissions are insertions and deletions of sites, also
called \textit{indels}. Stochastic models taking indels
into account have long been known~\cite{ThKiFe91,ThKiFe92}, but
they are not widely used in practice (or in theory) because
of their complexity. Instead, most practical algorithms
take a two-phase approach:
\begin{longlist}[(2)]
\item[(1)] \textit{Multiple sequence alignment.}
Site $t_i$ of sequence $\s_i$ and
site $t_j$ of sequence $\s_j$ are said to be \textit{homologous} if
they descend from the same site $t_0$ of a common ancestor $u$ (not
necessarily the MRCA)
\textit{only through substitutions}. In the \textit{multiple sequence
alignment} (MSA)
problem, we seek roughly to uncover the homology relation
between $\s_1,\ldots,\s_n$. Typically, the output is represented
by a matrix $\data$ of $n$ aligned sequences of equal length with
values in
$\{\A,\G,\C,\T,\gap\}$. Each column of the matrix corresponds to
homologous sites.
The state $\gap$ is called
a \textit{gap} and is used to account for insertions and deletions. For instance,
if sequence $\s_l$ does not have a site corresponding to $t_0$ in $u$ above,
then a gap is aligned with positions $t_i$ of $\s_i$ and
$t_j$ of $\s_j$ (which belong to the same column).

\item[(2)] \textit{Phylogenetic tree reconstruction.}
The matrix $\data$ is then cleaned up by removing all columns
containing gaps. Let $\data'$ be this new matrix. A standard
PTR algorithm is then applied to $\data'$. Note that substitutions
alone suffice to explain $\data'$.
\end{longlist}
Traditionally, most of the research on phylogenetic methods has focused
on the second phase.

In fact, current theoretical analyses of PTR assume that
the MSA problem has been solved \textit{perfectly}.
This has been a long-standing assumption in evolutionary biology.
But this simplification is increasingly being questioned in the
phylogenetic literature, where it has been argued that
alignment heuristics often create systematic
biases that affect analysis~\cite{LoytynojaGoldman08,WoSuHu08}.
Much recent empirical work has been devoted
to the proper joint estimation of alignments
and phylogenies \cite
{ThKiFe91,ThKiFe92,Metzler03,MiLuHo04,SuchardRedelings06,RivasEddy08,LoytynojaGoldman08,LRNLW09}.
Here we give the first analysis of an efficient, provably consistent
PTR algorithm in the presence of indels.
Our new algorithm suggests that a rough alignment
suffices for an accurate tree reconstruction (bypassing the
computationally difficult multiple alignment problem).\vspace*{10pt}

\textit{Theoretical properties of PTR.}\quad
In addition to computational efficiency, an important theoretical criterion
in designing a PTR algorithm is the so-called
\textit{sequence-length requirement} (SLR).
At a minimum, a reconstruction algorithm should be \textit{consistent},
that is, assuming a model of sequence evolution, the output
should be guaranteed to converge on the true tree as the sequence
length $k$ (the number of \textit{samples}) goes to $+\infty$ \cite
{Felsenstein78}.
Beyond consistency, the sequence-length requirement (or convergence
rate) of a
PTR algorithm is the sequence length required for guaranteed high-probability
reconstruction. The SLR is typically given as an asymptotic function of $n$,
the number of leaves of the tree. Of course, it also depends on the
substitution parameters.

A classical result due to Erd\H{o}s et al.~\cite{ErStSzWa99a} states that,
for general trees under the assumption that all branch lengths
are bounded by constants, the so-called short quartet method (SQM)
has $\poly(n)$-SLR. The SQM is a
particular PTR algorithm based on estimating evolutionary distances
between the leaf taxa, that is, the sum of the branch lengths between
species. Such algorithms are known as \textit{distance-based methods}.
The basic theoretical result behind distance-based methods is
the following: the collection of pairwise evolutionary distances
between all species forms a special metric on the leaves known
as an additive metric; under mild regularity assumptions,
such a metric \textit{characterizes} the underlying phylogeny
interpreted as an edge-weighted tree, that is, there is a one-to-one
correspondence between additive metrics and phylogenies;
moreover, the mapping between them can be computed efficiently
\cite{SempleSteel03}.\vspace*{10pt}

\textit{A new approach.}\quad
In the classical theoretical setting above where the MSA problem
is assumed perfectly solved
(we refer to this setting below as the ESSW framework),
the evolutionary distance between two species is measured using
the Hamming distance (or a state-dependent generalization)
between their respective sequences.
It can be shown that after a proper correction for multiple substitutions
(which depends on the model used) the expectation of the quantity
obtained does satisfy the additive metric property and can therefore
serve as the basis for a distance-based PTR algorithm.

Moving beyond the ESSW framework, it is tempting to account for
indels by simply using edit distance instead of the Hamming distance.
Recall that the \textit{edit distance} or
\textit{Levenshtein distance} between two strings is given by the
minimum number
of operations needed to transform one string into the other,
where an operation is an insertion, deletion or substitution of a
single character.
However, no analytical expression is known for the expectation of edit
distance under
standard indel models and computing such an expression appears
difficult (if
at all possible). An alternative idea is to compute the \textit{maximum
likelihood
estimator} for the time elapsed between two species given their
sequences. But this involves solving a nonconvex
optimization problem and the likelihood is only known to be efficiently
computable under a rather unrealistic assumption known as
reversibility~\cite{ThKiFe91}
(see below).

We use a different approach.
We divide the sequences into quantile blocks (the first $x \%$, the
second $x \%$, etc.).
We show that by appropriately choosing $x$ above
we can make sure that the blocks in different
sequences essentially ``match'' each other, that is,
they are made of mostly homologous sites.
We then compare the state frequencies in matching
blocks and build an additive metric out of this
statistic. As we show below, this is in fact
a natural generalization of the Hamming estimator
of the ESSW framework. However, unlike the Hamming distance
which can easily be analyzed through standard
concentration inequalities, proving rigorously
that our approach works involves several new
technical difficulties.
Our analysis relies on a branching process
analysis of the site displacements.
We give a quick proof sketch
after the formal statement of our results in Section~\ref{secresult}.

The results described here were first
announced without
proof in the special case of ultrametric trees under the CFN model with
inverse logarithmic indel rates~\cite{DaskalakisRoch10}. Here we give
full proofs
of stronger results,
including extensions to bounded-rate trees under
GTR models.\vspace*{8pt}

\textit{Related work.}\quad
For more background on models of molecular evolution and phylogenetics,
see, for example,~\cite{GraurLi99,SempleSteel03,Felsenstein04}.
Following the seminal results of~\cite{ErStSzWa99a}, there has been much
work on sequence-length requirement, including \cite
{Atteson99,ErStSzWa99b,HuNeWa99,SteelSzekely99,CsurosKao01,Csuros02,SteelSzekely02,KiZhZh03,MosselRoch06,DaMoRo06,LaceyChang06,DHJMMR06,Mossel07,GrMoSn08,Roch08,DaMoRo08a}.

The multiple sequence alignment problem as
a combinatorial optimization problem (finding the best
alignment under a given pairwise scoring function)
is known to be NP-hard~\cite{WangJiang94,Elias06}.
Most heuristics used in practice, such as
\mbox{CLUSTAL}~\cite{HigginsSharp88},
MAFFT~\cite{KaMiKuMi02}
and
MUSCLE~\cite{Edgar04},
use the idea of a guide tree, that is,
they first construct a very rough phylogenetic
tree from the data (using, e.g., edit distance
as a measure of evolutionary distance),
and then recursively construct
local alignments produced by ``aligning alignments.''

To our knowledge, little theoretical work has been
dedicated to the joint
estimation of alignments and phylogenies, with the exception
of Thatte~\cite{Thatte06} who gave consistency
results for the reversible case in the limit
where the deletion-to-insertion ratio tends to
1. However, no sequence-length requirement is obtained
in~\cite{Thatte06}. In recent related work,
the problem
of reconstructing ancestral sequences in the presence of
indels was considered~\cite{AnDaHaRo10,AnBrHa11}.\vspace*{-2pt}

\subsection{Model of sequence evolution}\vspace*{-2pt}

\subsubsection*{Phylogeny}
A \textit{phylogeny} is represented by a binary tree $T=(V,E)$,
whose leaves $L \subset V$ correspond
to extant species, and whose bifurcations
denote evolutionary events whereby two new species are generated
from an ancestor. The root of the phylogeny, denoted by $r(T)$,
represents the common ancestor of all the species in the phylogeny,
and we assume that all edges of $T$ are directed away from $r(T)$;
so, if $e=(u,v)$ is a branch of the phylogeny, $u$ is the \textit{parent}
of $v$ and $v$ is the \textit{child} of $u$. Moreover, if $v'$ is in
the subtree of $T$ rooted at $u$, we call $v'$ a \textit{descendant}
of $u$ and $u$ an \textit{ancestor} of $v'$.\vadjust{\goodbreak}

Along each branch of the phylogeny, the genetic material of the parent
species is subject to modifications that produce the genetic material
of its child species. A~common biological assumption is that the genetic
material of each species $u$ can be represented by a binary sequence
$\s_u= (\s_u^1,\ldots,\s_u^{K_u})$ of length $K_u$ over a finite
alphabet---for ease of presentation, we work with a binary alphabet $\{
0,1\}$ (but see Section~\ref{secextensions} for extensions to richer
alphabets)---and that the changes to which
$\s_u$ is subjected along the branch $e=(u,v)$ are described by a
Markov process.
In particular, the Markov property implies that, given the sequence $\s
_u$ of $u$,
the sequence $\s_v$ is independent of the sequences of the species outside
the subtree of $T$ rooted at $u$.

A simplifying assumption commonly used in phylogenetics is that all
species have sequences of the same length and, moreover, that every
\textit{site}, that is, every coordinate, in their sequences evolves independently
from every other site. In particular, it is assumed that, along each
branch $e=(u,v)$ of the phylogeny, every site $\s_u^j$ of the sequence
$\s_u$ is flipped with probability $p_e$ to the value $1-\s_u^j$
independently from the other sites. This model is known as the
Cavender--Farris--Neyman (CFN) model. A simple generalization
to $\{\A,\G,\C,\T\}$ is known as the Jukes--Cantor (JC) model
(see, e.g.,~\cite{Felsenstein04}).

\subsubsection*{Accouting for indels}
In this paper, we consider a more general evolutionary process that
accounts for the possibility of insertions and deletions.
Our model is similar to the original TKF91 model~\cite{ThKiFe91},
except that we do not enforce reversibility.
In our model, every edge $e=(u,v)$ of the phylogeny is characterized
by a quadruple of parameters $(t_e; \mutate_e, \death_e, \birth_e)$,
where $t_e$ is the evolutionary time between the species $u$ and $v$,
and $\mutate_e$, $\death_e$ and $\birth_e$ are, respectively, the
\textit{substitution}, \textit{deletion} and \textit{insertion} rates.
The Markov process by which the sequence at $v$ is obtained from the sequence
at $u$ is defined below (see, e.g.,~\cite{KarlinTaylor81} for
background on continuous-time
Markov processes).
%
\begin{definition}[(Evolutionary process on a branch)]
\label{defevolutionaryprocessbranch}
Given an edge $e=(u,v)$, with parameters $(t_e; \mutate_e, \death_e,
\birth_e)$,
the sequence $\s_v$ at $v$ is obtained from the sequence $\s_u$ at
$u$ according
to the following Markov process:
\begin{enumerate}[(2)]
\item[(1)] Intialize $\s_v: = \s_u$, $K_v:= K_u$ and $t_{\ell}
:=t_e$ (where $t_{\ell}$ is the remaining time on the edge $e$).

\item[(2)] While $t_{\ell} >0$:
\begin{itemize}
\item (\textit{Timing of next event})
let $I_0$, $I_1, \ldots, I_{K_v}$ be exponential random variables
with rate $\birth_e$, $D_1, \ldots, D_{K_v}$ exponential random variables
with rate $\death_e$ and $M_1, \ldots, M_{K_v}$ exponential
random variables
with rate $\mutate_e$; suppose that these random variables are
mutually independent
and let $\cal{T}$ be their minimum;

\item if ${\cal T} > t_{\ell}$, the process ends
at $t_\ell$; otherwise:
\begin{itemize}
\item(\textit{Insertion}) if $I_j = {\cal T}$, insert a new
site whose value is chosen uniformly at random from $\{0,1\}$\vadjust{\goodbreak} between
the sites
$\s^j_{v}$ and $\s^{j+1}_v$ of $\s_v$;
\item (\textit{Deletion}) if
$D_j = {\cal T}$, delete the site $\s_v^j$ from $\s_v$;
\item (\textit{Substitution}) and if
$M_j = {\cal T}$, replace $\s_v^j$ by $1-\s_v^j$;
\end{itemize}
(If $j=0$, then $\s^j_{v}$
is undefined and, if $j=K_v$, then $\s^{j+1}_{v}$ is undefined.)

\item (\textit{Remaining time}) update $\s_v$ according to these
changes, and update $K_v$ to reflect
the new sequence length; set the remaining time $t_{\ell}:=t_{\ell} -
{\cal T}$.
\end{itemize}
\end{enumerate}
\end{definition}

In words, the evolutionary process defined above assumes
that every site of the sequence $\s_u$ of the parent species is,
independently from the other sites, subjected to a sequence of
evolutionary events that flip its value; these events are distributed
according to a Poisson point process of intensity $\mutate_e$ in the
time interval $[0,t_e]$. However, the site may get deleted and
therefore not be inherited by the sequence of the node $v$;
this is determined by whether an exponential random variable
of rate $\death_e$ is smaller than $t_e$. While each site of
the parental sequence $\s_u$ is subjected to this process,
new sites are introduced in the space between existing sites
at rate~$\birth_e$, and each of these sites follows a similar
process for the remaining time. In essence,
insertion and deletion events are governed by an
independent branching process for each ancestral site.
Note further that the order of the sites, as described above, also
plays a role.
%
\begin{remark}
Unlike~\cite{ThKiFe91}, we do not use an ``immortal link''
and we do not assume that the length
process is at stationarity.
Our techniques can also be applied to the $\mathit{TKF}91$ model without much
modifications. We leave the details to the reader.
\end{remark}

Given the evolutionary process on a branch of the phylogeny,
the evolutionary process on the whole phylogeny is defined as follows.
%
\begin{definition}[(Evolutionary process)]\label{defevolutionaryprocess}
Suppose that every site of the sequence $\s_{r(T)}$ at the root of
the phylogeny is chosen to be $0$ or $1$ uniformly at random.
Recursively, if $\s_u$ is the sequence at node $u$ and $e=(u,v)$
is an edge of the phylogeny, the sequence $\s_v$ at node $v$
is obtained from the sequence $\s_u$ by an application of the
evolutionary process on a branch described by
Definition~\ref{defevolutionaryprocessbranch}.
\end{definition}

For ease of exposition, we first present our proof in the special case
where the
substitution, insertion and deletion rates are the same
on all edges of the phylogeny.
%
\begin{definition}[(Ultrametric assumption)]
Under the ultrametric assumption, the
leaves of the phylogeny are contemporaneous,
that is, there exists $H$ such that for each
$u \in L$ the sum of evolutionary times $t_e$
on the branches between $u$ and the root
is $H$.
\end{definition}
%
\begin{definition}[(Molecular clock assumption)] \label{defmolecularclock}
Under the molecular clock assumption,
we assume that the ultrametric assumption holds.
Moreover, there exist $\mutate$, $\death$ and
$\birth$ such that $\mutate_e = \mutate$,
$\death_e = \death$ and
$\birth_e=\birth$, for all $e \in E$.\looseness=-1
\end{definition}

We discuss a more general case in Section~\ref{secextensions}.

\subsubsection*{Notation}
In the sequel, we label the leaves of the phylogeny with the
positive integers $1$, $2, \ldots, n$, so that
$L = \{1,\ldots,n\}$, and the root $r(T)$ of the phylogeny with $0$.

\subsection{Main result}\label{secresult}

\subsubsection*{Statement of results}
We begin with a consistency result. Here
we consider a completely general phylogeny, that is,
neither the ultrametric nor the molecular clock
assumptions need hold.
%
\begin{theorem}[(Consistency: finite case)]\label{thmconsistency}
Assume that $0< t_e, \mutate_e,
\birth_e, \death_e < +\infty$, for all $e \in E$. 
Then there exists a procedure returning the correct tree
from the sequences at the leaves,
with probability of failure approaching $0$ as the sequence length at
the root of the tree goes to $+\infty$.
\end{theorem}

Our main result is the following. For simplicity
we first work under the symmetric two-state case and assume that
the molecular clock assumption holds.
%
\begin{theorem}[(Main result: two-state, molecular clock case)]\label{thmmain}
Consider the two-state model under the molecular clock assumption.
Assume further that
there exist constants
\[
0 < f,\qquad g <+\infty,
\]
independent of $n$,
such that
\[
f < t_e < g\qquad \forall e \in E.
\]
Moreover,
assume that
\[
\mutate_e = \mutate,\qquad
\birth_e= \birth,\qquad
\death_e = \death\qquad
\forall e\in E,
\]
where $\mutate$, $\birth$ and $\death$ are
bounded between constants (independent of $n$)
$0 < \underbar{\mutate} < \bar{\mutate} < +\infty$,
$0 = \underbar{\birth} < \bar{\birth} < +\infty$
and
$0 = \underbar{\death} < \bar{\death} < +\infty$,
respectively.
Under the assumptions above,
for all $\beta' > 0$ there exists
$\beta'' > 0$ such that
there exists a polynomial-time algorithm
solving the phylogenetic reconstruction problem
(i.e., returning the correct tree)
with probability of failure $n^{-\beta'}$,
if the root sequence has length
$k_r \geq n^{\beta''}$.\setcounter{footnote}{2}\footnote{In~\cite{DaskalakisRoch10}, a
preliminary version
of this result
was announced without proof, with
the much stronger assumption that
$\bar{\birth}, \bar{\death} = O(1/\log n)$,
that is, that the indel rates are negligible.
Here we show that this assumption
can be relaxed (at the cost of longer sequences).}\vadjust{\goodbreak}
\end{theorem}
%
\begin{remark}[(Branch lengths)]
Our assumption that all branch lengths $t_e$, $e \in E$, satisfy $f<
t_e < g$
is standard in the sequence-length requirement literature following the seminal
work of~\cite{ErStSzWa99a}.
\end{remark}
%

\subsubsection*{Extensions}
In Section~\ref{secextensions} we derive the following
extension. Let $Q$ be a reversible $4\times4$ rate
matrix with stationary
distribution $\pi$. (Larger alphabets are also possible.)
The GTR sequence evolution process is identical
to the one described in Definition~\ref{defevolutionaryprocessbranch}
except that the substitution process
is a continuous-time Markov process with rate matrix
$\mutate_e Q$.
%
\begin{theorem}[(Main result: GTR, bounded-rates case)]\label{thmmain2}
Consider the GTR model with rate matrix $Q$
under the
ultrametric assumption
(but not necessarily the molecular clock assumption).
Assume further that there exist constants
\[
0 < f,\qquad g,
\underbar{\mutate}, \bar{\mutate},
\underbar{\birth}, \bar{\birth},
\underbar{\death}, \bar{\death},
<+\infty,
\]
independent of $n$,
such that
\[
f< t_e < g,\qquad
\underbar{\mutate} < \mutate_e
< \bar{\mutate}\qquad
\forall e \in E.
\]
Moreover,
assume that
\[
\birth_e= \birth,\qquad
\death_e = \death\qquad
\forall e\in E,
\]
where $\birth$ and $\death$ are
bounded between constants (independent of $n$)
$0 = \underbar{\birth} < \bar{\birth} < +\infty$
and
$0 = \underbar{\death} < \bar{\death} < +\infty$,
respectively. We refer to the conditions above as
the \textit{bounded-rates assumption}.
Under the assumptions above,
for all $\beta' > 0$ there exists
$\beta'' > 0$ such that
there exists a polynomial-time algorithm
solving the phylogenetic reconstruction problem
(i.e., returning the correct tree)
with probability of failure $n^{-\beta'}$,
if the root sequence has length
$k_r \geq n^{\beta''}$.
\end{theorem}

\subsubsection*{Proof sketch}
Consider the two-state, molecular clock case.
As we noted before, unlike the classical setting where the
Hamming distance can be analyzed through standard
concentration inequalities, proving rigorously
that our approach works involves several new
technical difficulties.
The proof goes through the following
steps:
\begin{longlist}[(4)]
\item[(1)] \textit{Expectations.} We first compute expectations
of block statistics, which involve analyzing a continuous-time
Markov process. We use these calculations to define an appropriate
additive metric based on correlations between blocks.

\item[(2)] \textit{Sequence length and site displacements.} We give bounds
on how much
sequence lengths vary across the tree
through a moment-generating function argument.
Using our bounds on the sequence length
process, we bound the worst-case displacements of the sites. Namely, we show
that, under our assumptions, all sites move by at most $O(\sqrt{k \log k})$.

\item[(3)] \textit{Sequence partitioning.}
We divide each sequence in blocks of size rough\-ly $k^{\zeta}$ for
$\zeta> 1/2$, where $k$ is the sequence length at the root. From our bounds
on site displacements, it follows that the blocks roughly
match across different sequences. In particular, we bound
the number of homologous sites between matching blocks with high
probability and
show that the expected correlation between these blocks
is approximately correct.

\item[(4)] \textit{Concentration.} Finally, we show that our estimates are
concentrated.
The concentration argument proceeds by conditioning on the indel process
satisfying the high-probability conditions in the previous points.
\end{longlist}
The crux of our result is the proper estimation
of an additive metric.
With such an estimation procedure in hand,
we can use a standard distance-based approach
to recover the phylogeny.

\subsubsection*{Organization}
The rest of the paper is organized as follows.
The evolutionary distance forming the basis of our approach
is presented in Section~\ref{secmetric}.
We describe our full distance estimator in Section~\ref{secestimation}
and prove its concentration in the same section.
Extensions are described in Section~\ref{secextensions}.

\section{Evolutionary distances}\label{secmetric}

Consider the two-state, molecular clock case.
In this section, we show how to define an appropriate
notion of ``evolutionary distance'' between two species.
Although such distances have been widely used in prior phylogenetic
work and have been defined for a variety of models \cite
{SempleSteel03,Felsenstein04},
to our knowledge our definition is the first that applies
to models with indels.
We begin by reviewing the standard definition in the indel-free case
and then adapt it to the presence of indels. Our estimation procedure
is discussed in Section~\ref{secestimation}.

\subsection{The classical indel-free case}

Suppose first that $\birth= \death= 0$,
that is,
\textit{there is no indel}.
In that case, the sequence length remains fixed at
$k$ and the alignment problem is trivial.
Underlying all distance-based approaches is the following
basic definition.
%
\begin{definition}[(Additive metric)]\label{defadditive}
A phylogeny is naturally equipped with a so-called additive metric
on the leaves
$\dist\dvtx L\times L \to(0,+\infty)$ defined as
\[
\forall a,b \in L\qquad \dist(a,b) = \sum_{e\in\path_T(a,b)} \weight_e,
\]
where $\path_T(a,b)$ is the set of edges on the path between
$a$ and $b$ in $T$ and where $\weight_e$ is a nonnegative function of the
parameters on $e$ (in our case, $t_e$, $\mutate_e$, $\birth_e$
and~$\death_e$). For instance, a common choice for $\weight_e$ would be
$\weight_e = \mutate_e t_e$ in which case
$\dist(a,b)$ is
the expected number of substitutions per site between
$a$ and $b$.
Often $\dist(a,b)$ is referred to as the ``evolutionary distance''
between species $a$ and $b$.
Additive metrics\vadjust{\goodbreak} are characterized by the following four-point condition:
for all \mbox{$a,b,c,d \in L$},
\[
\dist(a,b) + \dist(c,d) \leq\max\{\dist(a,c) + \dist(b,d), \dist
(a,d) + \dist(b,c)\}.
\]
Moreover, assuming $\weight_e > 0$ for all $e\in E$,
it is well known that there exists
a one-to-one correspondence between $\dist$ and $T$ as a weigthed tree
with edge weights $\{\weight_e\}_{e\in E}$. We will discuss algorithms
for constructing $T$ from $\dist$ in Section~\ref{secalgo}.
For more background on tree-based metrics, see~\cite{SempleSteel03}.
\end{definition}

Definition~\ref{defadditive} implies that phylogenies can be
reconstructed by computing $\dist(a,b)$ for all pairs of leaves $a,b\in
L$. Assume we seek to estimate the evolutionary distance between
species $a$ and $b$ using their respective sequences. In a first
attempt, one might try the (normalized) Hamming distance between $\s_a
= (\s_a^1,\ldots,\s_a^k)$ and $\s_b = (\s_b^1,\ldots,\s_b^k)$. However,
the expected Hamming distance (in other words, the probability of
disagreement between a site of $a$ and~$b$) does not form an additive
metric as defined in Definition~\ref{defadditive}. Instead, it is well
known that an appropriate estimator is obtained by ``correcting''
the\vspace*{1pt} Hamming distance for ``multiple'' substitutions.
Denoting by $\hamming(\s_a,\s_b)$ the Hamming distance between $\s_a$
and $\s_b$, a Markov chain calculation shows that $\dist(a,b) =
-\frac{1}{2}\log(1 - 2\expec[\hamming(\s_a,\s_b)])$, with the choice
$\weight_e = \mutate_e t_e$ (see, e.g.,~\cite{Felsenstein04}). In a
distance-based reconstruction procedure, one first estimates $\dist$
with
%
\begin{equation}\label{eqhamming}
\edist(a,b) = -\tfrac{1}{2}\log\bigl(1 - 2\hamming(\s_a,\s_b)\bigr)
\end{equation}
and then applies one of the algorithms discussed in Section \ref
{secalgo} below. The sequence-length\vspace*{1pt} requirement of such a
method can be derived by using concentration results for $\hamming$
\cite {ErStSzWa99a,Atteson99}.

\subsection{Taking indels into account}

To simplify the presentation, we assume
throughout that $\birth\neq\death$.
The case $\birth=\death$ follows from the
same argument.

In the presence of indels, the estimator (\ref{eqhamming}) based on
the Hamming
distance is difficult to apply. One has to first align
the sequences, which
cannot be done perfectly
and causes biases as well as correlations
that are hard to analyze.
Alternatively, one could try a different string distance such as edit
distance. However, computing the expectation of edit distance under indel
models appears difficult.

We use a different approach involving correlations
between state frequencies. We will eventually
apply the estimator to large sub-blocks of the sequences (see
Section~\ref{secestimation}),
but we first describe it for the full sequence for clarity.
For a node $u$, let
$K_u$ be the (random) length of the sequence at $u$ and
$Z_u$, the number of $0$'s in the sequence at $u$.
Then, our distance estimator is\looseness=-1
\[
\edist(a,b) = \bigl(Z_a - \tfrac{1}{2} K_a\bigr)\bigl(Z_b - \tfrac
{1}{2} K_b\bigr).
\]\looseness=0
We now analyze the expectation of this quantity.
For $u \in V$, we let
\[
\deviate_u = Z_u - \tfrac{1}{2} K_u\vadjust{\goodbreak}
\]
be the deviation of $Z_u$ from its expected value (conditioned
on the sequence length).

\subsubsection*{Single channel}
Suppose $T$ is made of a single edge from the
root $r$ to a leaf $a$ with parameters $t, \mutate, \birth, \death$.
Assume first that the original sequence length is $k_r = 1$.
Let $K_a$ be the length of the sequence at $a$.
Then $K_a$ is a continuous-time
branching process and, by Markov chain
calculations (\cite{AthreyaNey72}, Section III.5),
its moment-generating function is
%
\begin{equation}\label{eqmomentgeneratingfunction}
F(s,t) \equiv\expec[s^{K_a}]
=
\frac{\death(s-1) - e^{(\death-\birth)t}(\birth s - \death)}
{\birth(s-1) - e^{(\death-\birth)t}(\birth s - \death)}.
\end{equation}
By differentiating $F(s,t)$ we derive
%
\begin{equation}\label{eqgrowth}
\expec[K_a]
=
e^{-(\death-\birth)t}
\end{equation}
and
%
\begin{equation}\label{eqgrowth2}
\var[K_a]
=
\frac{\death+\birth}{\death-\birth}\bigl[e^{-(\death-\birth)t} -
e^{-2(\death-\birth)t}\bigr].
\end{equation}
Let $K^*_a$
be the number of ``new'' sites at $a$, that is, excluding the original
site if it survived. (We ignore the substitutions for the time being.)
The probability that the original site survives is $e^{-\death t}$.
Then,
\[
\expec[K^*_a]
= \expec[K_a - \ind\{\mbox{original site survives}\}]
= e^{-(\death-\birth)t} - e^{-\death t}
\]
by linearity of expectation.

We now take into account substitutions.
Assume that the original sequence length at $r$ is a random variable
$K_r$ and that the sequence at $r$ is i.i.d. uniform. Denote by $Z_r$
the number of
$0$'s at $r$.
The probability that a site in $r$, that is still surviving in $a$, has flipped
its value is
\begin{eqnarray*}
p
&=&
\prob[ \mbox{state flips odd number of times in time $t$} ]\\
&=& \sum_{j=0}^{+\infty} e^{-\mutate t} \frac{(\mutate
t)^{2j+1}}{(2j+1)!}\\
&=& e^{-\mutate t}\sinh\mutate t\\
&=& \frac{1 - e^{-2\mutate t}}{2}.
\end{eqnarray*}
Also, note that a new site created along the path between $r$ and $a$
has equal chance of being $0$ or $1$ \textit{at the end of the path}.
Then we have the following lemma.
%
\begin{lemma}[(Single channel: expected deviation)]\label{lem1}
The following holds:
\[
\expec[\deviate_a | K_r,Z_r] = e^{-(2\mutate+\death) t} \deviate_r.\vadjust{\goodbreak}
\]
\end{lemma}
\begin{pf}
We have
%
\begin{eqnarray}\label{eqsingle}
\expec[\deviate_a | K_r,Z_r]
&=& \expec\bigl[\bigl(Z_a - \tfrac{1}{2} K_a\bigr) | K_r,
Z_r\bigr]\nonumber\\
&=& Z_r e^{-\death t} (1-p)
+ (K_r - Z_r) e^{-\death t} p\nonumber\\
&&{} + K_r \bigl(e^{-(\death-\birth)t} - e^{-\death t}\bigr) \tfrac{1}{2}
- K_r e^{-(\death-\birth)t} \tfrac{1}{2}\\
&=& Z_r (1- 2p) e^{-\death t}
- \tfrac{1}{2} K_r (1-2p) e^{-\death t}\nonumber\\
&=& e^{-2\mutate t} e^{-\death t} \deviate_r,\nonumber
\end{eqnarray}
where on the first two lines:
\begin{longlist}[(3)]
\item[(1)]
the first term is the number of original
$0$'s surviving in state $0$;
\item[(2)]
the second term is the number of original
$1$'s surviving in state $0$;
\item[(3)]
the third term is the number of new sites surviving in state $0$ (where
recall that new sites
are uniformly chosen in $\{0, 1\}$);
\item[(4)]
the fourth term is half the sequence length at $a$ given
the length at~$r$.\qed
\end{longlist}
\noqed\end{pf}

\subsubsection*{Fork channel}
Consider now a ``fork'' tree, that is,
a root $r$ from which emanates a single edge $e_u = (r,u)$ which in turn
branches into two edges $e_a = (u,a)$ and $e_b = (u,b)$
(see Figure~\ref{figfork} below).
For $x=a,b,u$, we denote the parameters of edge $e_x$ by
$t_x, \birth_x, \death_x, \mutate_x$.
Our goal is to compute $\expec[\edist(a,b)]$ assuming that the
sequence length
at the root is $k_r$. We use (\ref{eqsingle}),
the Markov property and
the fact that $Z_u$ conditioned on $K_u$ is a binomial
with parameters $1/2$ and $K_u$. We get the following lemma.
%
\begin{lemma}[(Fork channel: expected distance)]\label{lem2}
The following holds:
\[
\expec[\edist(a,b)]
= e^{-(2\mutate_a+\death_a) t_a}
e^{-(2\mutate_b+\death_b) t_b}
e^{-(\death_u - \birth_u)t_u} \frac{k_r}{4}.
\]
\end{lemma}
\begin{pf}
We have
\begin{eqnarray*}
\expec[\edist(a,b)]
&=& \expec[\deviate_a \deviate_b]\\
&=& \expec[\expec[\deviate_a \deviate_b |
K_u,Z_u]]\\
&=& \expec[\expec[\deviate_a | K_u,Z_u]\expec
[\deviate_b | K_u,Z_u]]\\
&=& e^{-2\mutate_a t_a} e^{-\death_a t_a} e^{-2\mutate_b t_b}
e^{-\death_b t_b} \expec[\deviate_u^2]\\
&=& e^{-2\mutate_a t_a} e^{-\death_a t_a} e^{-2\mutate_b t_b}
e^{-\death_b t_b} \expec[\expec[\deviate_u^2 |
K_u]]\\
&=& e^{-2\mutate_a t_a} e^{-\death_a t_a} e^{-2\mutate_b t_b}
e^{-\death_b t_b} \expec\biggl[\frac{K_u}{4}\biggr]\\
&=& e^{-2\mutate_a t_a} e^{-\death_a t_a}
e^{-2\mutate_b t_b} e^{-\death_b t_b}
\frac{e^{-(\death_u - \birth_u)t_u} k_r}{4},
\end{eqnarray*}
where we used (\ref{eqgrowth})
and Lemma~\ref{lem1}.
\end{pf}

\subsubsection*{Molecular clock}
We specialize the previous result to the molecular clock
assumption. That is, we assume, for $x=a,b,u$,
that $\birth_x = \birth$, $\death_x = \death$
and $\mutate_x = \mutate$. Note that by construction
$t_a = t_b$ (assuming species $a$ and $b$ are contemporary).
We denote $t = t_a$ and $\maxt= t_u + t_a$.
Denoting
$\kappa= \frac{k_r e^{-(\death- \birth)\maxt}}{4}$, we then get
the following lemma.
%
\begin{lemma}[(Molecular clock: expected distance)]\label{lem3}
The following holds:
\[
\expec[\edist(a,b)]
= e^{-(4\mutate+\death+\birth)t}
\kappa.
\]
\end{lemma}

Letting
\[
\beta= 4\mutate+\death+\birth,
\]
we get that
\[
-2 \log\expec[\kappa^{-1} \edist(a,b)] = 2\beta t,
\]
which is the evolutionary distance between $a$ and $b$
with the choice $\weight_e = \beta t_e$.
Therefore, we define the following estimator:
\[
\edist^*(a,b)
= -2 \log\kappa^{-1} \edist(a,b).
\]

\section{Distance computation}\label{secestimation}

We now show how to estimate the evolutionary distance between
two species by decomposing the sequences into large blocks
which serve as roughly independent samples. We use the
following notation:
$M_t = e^{-(\death-\birth)t}$,
$D_t = e^{-\death t}$,
$\netrate= \death- \birth$,
$\totalrate= \death+ \birth$
and
$\Gamma_t = \netrate^{-1} \birth(1 - M_t)$.

We show in Section~\ref{secalgo}
that the time elapsed between the root
and the leaves is bounded by
$\frac{g^2}{f} \log_2 n$.
Hence, under our assumptions
%
\begin{eqnarray}\label{eqboundm}
\Upsilon_n^{-1}
&\equiv&
e^{-(\bar{\death} +\bar{\birth}) ({g^2}/{f}) \log_2 n}
\leq
e^{-(\bar{\death} {g^2}/{f}) \log_2 n}\nonumber\\[-8pt]\\[-8pt]
&\leq& M_t
\leq
e^{(\bar{\birth} {g^2}/{f}) \log_2 n}
\leq
e^{(\bar{\birth}+\bar{\death}) ({g^2}/{f}) \log_2 n}
\equiv\Upsilon_n,\nonumber
\\
%
\label{eqboundd}
\Upsilon_n^{-1}
&\leq&
e^{-(\bar{\death} {g^2}/{f}) \log_2 n}
\leq D_t
\leq1
\end{eqnarray}
and
%
\begin{eqnarray}\label{eqboundgamma}
0
&\leq&\Gamma_t
= \birth t \frac{1 - e^{-(\death- \birth)t}}{(\death- \birth)t}
\leq
\bar{\birth} \frac{g^2}{f} \log_2 n
\frac{e^{(\bar{\birth} {g^2}/{f}) \log_2 n} - 1 }{(\bar{\birth}
{g^2}/{f}) \log_2 n} \nonumber\\[-8pt]\\[-8pt]
&=& e^{(\bar{\birth} {g^2}/{f}) \log_2 n}
- 1
\leq\Upsilon_n,\nonumber
\end{eqnarray}
where we used that the function $x^{-1}(1 - e^{-x})$
is nonnegative and decreasing since its derivative is
\[
\frac{x e^{-x} - (1 - e^{-x})}{x^2}
= e^{-x} \frac{(1+x) - e^x}{x^2}
\leq0,\qquad x\neq0.
\]
Note that the bounds above are polynomials in $n$ with exponents
depending only on $f$, $g$, $\bar{\birth}$ and $\bar{\death}$. In
particular,\vspace*{2pt} we\vadjust{\goodbreak} will ultimately take sequence lengths $k_r$
of the form $n^{\beta''}$ with $\beta''$ chosen much larger than the
exponent in $\Upsilon_n$. We call polynomials in $n$ (such as
$\Upsilon_n$) which have an exponent not depending on $\beta''$,
\textit{small polynomials}. As a result, the following notation will be
useful. For a function $W(k_r)$ of $k_r$, we use $\scal_n(W(k_r))$ to
denote a function smaller or equal to $W(k_r)$ \textit{up to a small
polynomial factor}. (The latter will be used similarly to the big-O
notation.)

Recall the following standard concentration inequalities (see,
e.g.,~\cite{MotwaniRaghavan95}).

\begin{lemma}[(Chernoff bounds)]\label{lemmachernoff}
Let $Z_1,\ldots,Z_m$ be independent $\{0,1\}$-random
variables
such that, for $1 \leq i \leq m$,
$\prob[Z_i = 1] = p_i$
where $0 < p_i < 1$. Then, for
$Z = \sum_{i=1}^m Z_i$, $M = \expec[Z] = \sum_{i=1}^m p_i$,
$0 < \delta_- \leq1$ and $0 < \delta_+ \leq U$,
\[
\prob[Z < (1-\delta_-) M] < e^{-M \delta_-^2/2}
\]
and
\[
\prob[Z > (1+\delta_+) M] < e^{-c(U) M \delta_-^2},
\]
where $c(U) = [(1+U) \ln(1+ U) - U]/U^2$.
\end{lemma}

\subsection{Concentration of the indel process} \label{seclengthconcentration}

\subsubsection*{Sequence length}
We first show that the sequence length is concentrated.
Let $T$ be single channel consisting of edge $e=(r,a)$. Let $k_r$ be
the length at~$r$.

\begin{lemma}[(Single channel: large deviations of sequence
length)]\label{lemmalengthconcentration}
For all $\gamma>0$ and
${\widehat{k}_r} \ge k_r = n^{\beta'''}$ with
$\beta''' > 0$ large enough,
with probability at least $1- {\widehat{k}_r}^{-\gamma}$,
\[
K_a = k_r M_t
\pm\scal_n\bigl(\sqrt{\widehat{k}_r \log{\widehat{k}_r}}\bigr),
\]
where the small polynomial factor in
$\scal_n(\sqrt{\widehat{k}_r \log{\widehat{k}_r}})$
depends on $\gamma$ as well.
\end{lemma}
%
\begin{remark}\label{remblocks}
Although we stated Lemma~\ref{lemmalengthconcentration} for the
full sequence,
it will also be needed for ``half-sequences'' and ``blocks.''
In particular, we use the previous lemma to track the position
of sites. In that context, one should think of $k_r$ as the
position of a site in $r$ and $K_a$ as its position in
$a$. Then we can use ${\widehat{k}_r}$ for the full sequence
length at $r$
(see Section~\ref{secpartition}).
\end{remark}
\begin{pf*}{Proof of Lemma~\ref{lemmalengthconcentration}}
We think of $K_a$ as
\[
K_a = \sum_{i=1}^{k_r} K_{a,i},
\]
where $K_{a,i}$ is the number of sites generated by a single
site of the sequence at $r$. Intuitively, $K_{a,i}$ is
the number of sites that were inserted between\vadjust{\goodbreak} the sites
$i$ and $i+1$ of the sequence at $r$, plus the site at
position $i$ itself, if it survived. Clearly the variables
$\{K_{a,i}\}_i$ are mutually independent.

Using \eqref{eqgrowth}
we obtain that
\[
\expec[K_a] = k_r M_t.
\]
%
For $\eps> 0$, by Markov's inequality, we have
%
\begin{equation}\label{eqmarkov}
\prob[K_a \geq k_r M_t + k_r \eps]
\leq s^{-k_r(M_t + \eps)} \expec[s^{K_{a}}]
= \bigl(s^{-(M_t + \eps)} \expec[s^{K_{a,1}}]\bigr)^{k_r}.
\end{equation}
We take $s = 1 + C \eps$ for $C > 0$ to be determined.

We have
\begin{eqnarray*}
\expec[s^{K_{a,1}}]
&=& \frac{\death(s-1) - e^{(\death-\birth)t}(\birth s - \death)}
{\birth(s-1) - e^{(\death-\birth)t}(\birth s - \death)}
= \frac{(\death- \birth M_t^{-1})C \eps
+ \netrate M_t^{-1}}
{\birth(1 - M_t^{-1})C \eps
+ \netrate M_t^{-1}}\\
&=& \frac{\netrate^{-1}(\death M_t - \birth)C \eps
+ 1}
{\netrate^{-1}\birth(M_t - 1)C \eps
+ 1}
= \frac{1 - (\birth^{-1} \death\Gamma_t - 1) C \eps}{1 - \Gamma
_t C \eps}\\
&=& [1 - (\birth^{-1} \death\Gamma_t - 1) C \eps]
\sum_{\iota= 0}^{+\infty} [\Gamma_t C \eps]^\iota,
\end{eqnarray*}
whenever $\Gamma_t C \eps< 1$.
Hence, if $\Upsilon_n C \eps< 1$ is bounded away from
$1$ (independently of $n$),
we have, using (\ref{eqboundgamma}),
\begin{eqnarray*}
\expec[s^{K_{a,1}}]
&=& [1 - (\birth^{-1} \death\Gamma_t - 1) C \eps]
[1 + \Gamma_t C \eps+ (\Gamma_t C \eps)^2
+ O((\Upsilon_n C \eps)^3)]\\
&=& 1 + M_t (C\eps) + M_t \Gamma_t (C\eps)^2
+ O((\Upsilon_n C \eps)^3).
\end{eqnarray*}
Moreover, using the binomial series and (\ref{eqboundm}),
and assuming $C\eps< 1$
\begin{eqnarray*}
s^{-(M_t + \eps)}
&=& \sum_{\iota=0}^{+\infty}
\frac{(-M_t - \eps) (-M_t - \eps-1) \cdots(-M_t - \eps- \iota+
1)}{\iota!}[C \eps]^\iota\\
&\leq& 1 - (M_t + \eps) (C\eps) + \frac{(M_t+\eps)(M_t + \eps+
1)}{2} (C\eps)^2\\
&&{} +
\sum_{\iota=3}^{+\infty}
(M_t + \eps+ 1)^\iota[C \eps]^\iota\\
&=& 1 - (M_t + \eps) (C\eps) + \frac{(M_t+\eps)(M_t + \eps+ 1)}{2}
(C\eps)^2\\
&&{}+ O((\Upsilon_n C \eps)^3),
\end{eqnarray*}
whenever $\eps$ is small and $\Upsilon_n C \eps< 1$ is bounded away
from $1$ (independently
from~$n$).
Therefore,
\begin{eqnarray*}
s^{-(M_t + \eps)} \expec[s^{K_{a,1}}]
&=& 1 - \eps(C\eps)+ M_t\Gamma_t(C\eps)^2
+ \frac{(M_t+\eps)(M_t + \eps+ 1)}{2} (C\eps)^2\\
&&{}
- (M_t + \eps)M_t (C\eps)^2
+ O((\Upsilon_n C \eps)^3).
\end{eqnarray*}
Note that the second term on the right-hand side depends
on $C$ whereas the remaining terms depend on $C^2$.
Taking $C = \Upsilon_n^{-2} C_0(\gamma)$
with $C_0(\gamma) > 0$ small enough and
$c = c_0(\gamma) > 0$
large enough, using (\ref{eqmarkov})
with the choice
\[
\eps= c\sqrt{\frac{\Upsilon_n^2 \log{\widehat{k}_r}}{k_r}},
\]
we get that
\begin{eqnarray*}
\prob\bigl[K_a \geq k_r M_t + c \sqrt{\Upsilon_n^2 {\widehat{k}_r}
\log{\widehat{k}_r}}\bigr]
&\leq&
\prob[K_a \geq k_r M_t + k_r \eps]\\[-2pt]
&\leq& \biggl(1 - \frac{O(\log\widehat{k}_r)}{k_r}\biggr)^{k_r}\\[-2pt]
&\leq& {\widehat{k}_r}^{-\gamma}.
\end{eqnarray*}
Note that our choice of $\eps$ satisfies
$\Upsilon_n C \eps< 1$ for $k_r$ a large enough
polynomial of~$n$ (compared to the small
polynomial $\Upsilon_n$).

A similar inequality holds for the other direction.
\end{pf*}

\subsubsection*{Correlated sites} 
\begin{figure}

\includegraphics{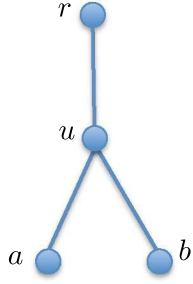}

\caption{The fork channel.}\label{figfork}
\end{figure}

Now let $T$ be the fork channel
consisting of nodes $r$, $u$, $a$ and $b$ as in Figure~\ref{figfork}.
%
Assume that $a$ and $b$ are contemporary, call $t$ the time separating
them from $u$ and denote by $\js{ab}$ the number of sites
in $a$ and $b$ that are jointly surviving from $u$. These are the sites that
produce correlation between the sequences at $a$ and $b$. All other sites
are essentially noise. We bound the large deviations of~$\js{ab}$.
%
\begin{lemma}[(Fork channel: large devations of jointly surviving
sites)]\label{lemsurviving} Condition on the sequence length at $u$
being $k_u$. Then,\vspace*{1pt} for all $\gamma>0$ and all
${\widehat{k}_u} \ge k_u = n^{\beta'''}$ with $\beta''' > 0$ large
enough, with conditional probability at least $1-
{\widehat{k}_u}^{-\gamma}$,
\[
\js{ab} = k_u D_t^2
\pm\scal_n\bigl(\sqrt{\widehat{k}_u \log{\widehat{k}_u}}\bigr),
\]
where the small polynomial factor in
$\scal_n(\sqrt{\widehat{k}_u \log{\widehat{k}_u}})$
depends on $\gamma$ as well.\vadjust{\goodbreak}
\end{lemma}
\begin{pf}
Each site in $u$ survives in $a$ with probability $D_t$.
The same holds for $b$ independently.

The result then follows from
Chernoff's bound. We have
\begin{eqnarray*}
\prob\bigl[
\js{ab} < k_u D_t^2 - c \sqrt{\Upsilon_n^2 \widehat{k}_u \log
\widehat{k}_u}
\bigr]
&\leq&
\prob\Biggl[
\js{ab} < k_u D_t^2 - c \sqrt{\frac{\Upsilon_n^2 \log\widehat
{k}_u}{k_u}} \cdot k_u D_t^2
\Biggr]\\
&\leq&
\exp(
- c^2 \Upsilon_n^2 D_t^2 \log\widehat{k}_u
)\\
&\leq& \widehat{k}_u^{-\gamma}
\end{eqnarray*}
for $c = c(\gamma) > 0$ large enough,
where we used (\ref{eqboundd}).

The other direction is similar.
\end{pf}

\subsection{Sequence partitioning}\label{secpartition}

From Lemma~\ref{lemmalengthconcentration}, it follows that the sites
of the root sequence (or of an internal sequence)
remain fairly close to their expected position at the leaves.
We take advantage
of this fact by dividing each sequence into blocks
of size asymptotically larger than the typical displacement
implied by Lemma~\ref{lemmalengthconcentration}.
As a result, matching blocks in different sequences share
a significant fraction of sites. Moreover, distinct blocks
are roughly independent. We estimate the evolutionary
distance between two leaves by comparing the
site frequencies in matching blocks. This requires
some care as we show next.

Consider the fork channel.
We seek to estimate the evolutionary distance $\edist(a,b)$
between $a$ and $b$ (normalized by the sequence length at $u$).

\subsubsection*{Partitioning the leaf sequences}
Let $k_0$ be some \textit{deterministic} length (to be determined), and
consider the first $k_0$ sites in the sequences $\s_a$ and $\s_b$ at
the nodes $a$ and $b$, respectively. If the sequence at $a$ or $b$ has
length smaller than $k_0$, we declare that our distance estimate
$\epartdist(a,b)$ (see below) is $+\infty$.

We divide the leaf sequences into $L$ blocks
of length $\ell$ where $\ell= \lceil k_0^{\zeta} \rceil$, for some
$\frac{1}{2} < \zeta<1$ to be determined later and
$L = \lfloor k_0 / \ell\rfloor$. We let $k_0'=\ell L$.
For all $i=1,\ldots, L$, we define the $i$th \textit{block} $\s_{a,i}$
of $a$
to be the subsequence of $\s_a$
ranging from position $(i-1) \ell+1$ to position $i \ell$.
We let $Z_{a,i}$ be the number of zeros inside $\s_{a,i}$ and
define the \textit{block deviations}
\[
\deviate_{a,i}=Z_{a,i} - {\ell\over2}
\]
for all $i=1,\ldots,L$, and similarly for the sequence at $b$.

Using the above notation we define our distance estimator next.
Assume that $L$ is even.
Otherwise, we can just drop the last block in the above partition.
Our estimator is
\[
\epartdist(a,b) = {2 \over L}
\sum_{j=0}^{L/2-1}\deviate_{a,2j+1} \deviate_{b,2j+1}.
\]
Notice that in our summation above we skipped every other block
in our sequence partition to avoid overlapping sites and hence,
decrease potential correlations between the terms in the
estimator.
In the rest of this section, we analyze
the properties of $\epartdist(a,b)$. To do this it is helpful
to consider the sequence at $u$ and the events that happened
in the channels defined by the edges $(u,a)$ and $(u,b)$.

\subsubsection*{Partitioning the ancestral sequence}
Let us choose $\ell_u$ to be the largest integer satisfying
%
\begin{equation}\label{eqellu}
{\ell}_u M_t \le\ell.
\end{equation}
Suppose that the sequence $\s_u$ at node $u$ is not shorter
than $k_u'=(L-1) \ell_u$, and define the $i$th \textit{ancestral block}
$\s_{u,i}$ of $u$ to be the
subsequence of $\s_{u}$ ranging from position
$(i-1) \ell_u +1$ to position $i \ell_u$, for all $i \le{L-1}$.
Given Lemma~\ref{lemmalengthconcentration}, the choice of $\ell_u$
in \eqref{eqellu} is such that the blocks of $u$ and the
corresponding blocks at $a$ and $b$ roughly align.

In order to use the expected evolutionary distance as computed in
Lem\-ma~\ref{lem3}, we define an ``interior'' ancestral block which
is guaranteed with high probability to remain entirely
``inside'' the corresponding leaf block.
Let $\delta_u= \lceil L + {1 \over M_t} \scal_n(\sqrt{k_u' \log
k_u'}) \rceil$,
where the small polynomial factor is the maximum of
those in the proofs of
Lemma~\ref{lemmalengthconcentration} and Lemma~\ref{lemsurviving}
for a given choice of $\gamma$.
[The $L = o(\sqrt{k_0})$ in $\delta_u$ is needed
only when \eqref{eqellu} is a strict inequality. See the proof of
Lemma~\ref{lemprobabilityofgoodevent1}
below.]
We define the $i$th (\textit{ancestral}) \textit{interior block} $\s_{u,i}'$
of $u$ to be the subsequence of $\s_{u,i}$ ranging from
position $(i-1) \ell_u + \delta_u$ of $\s_u$ to position
$i \ell_u - \delta_u$.
Notice that $\delta_u = \scal_n( \sqrt{k_0 \log k_0})$,
while $\ell_u = \scal_n(k_0^{\zeta})$.
Therefore, for $k_0 > k_0^*$,
where $k_0^*$ is sufficiently large,
$(i-1) \ell_u + \delta_u < i \ell_u - \delta_u$ so that the sequence
$\s_{u,i}'$ is well defined.

Also, for all $i=1,\ldots,L-1$, we define $x'_{a,i}$, $y'_{a,i}$
to be the position of the left-most (resp., right-most) site
in the sequence $\s_a$ descending from the site at position
$(i-1) \ell_u + \delta_u$ (resp., $i \ell_u - \delta_u$ of $\s_u$).
Similarly, we define $x'_{b,i}$ and $y'_{b,i}$.
Given this notation, we define the following ``good'' event
%
\begin{eqnarray}\label{eqgoodevent}
&&\good'_1 = \{
\forall i \le L-1\dvtx
(i-1) \ell< x'_{a,i},x'_{b,i} < (i-1) \ell+ 2M_t\delta_u,\nonumber
\nonumber\\[-8pt]\\[-8pt]
&&\hspace*{139.7pt} i \ell-
2M_t\delta_u < y'_{a,i},y'_{b,i} < i \ell
\}.\nonumber
\end{eqnarray}
Intuitively, when the event $\good'_1$ holds, all
surviving descendants of the
interior block $\s_{u,i}'$ are located inside the blocks
$\s_{a,i}$ and $\s_{b,i,}$, respectively
(and the blocks remain large enough).

To argue about block independence, we also define the \textit{exterior block}
$\s_{u,i}''$
of $u$ to be the subsequence of $\s_{u,i}$ ranging from
position $(i-1) \ell_u - \delta_u$ of $\s_u$ to position
$i \ell_u + \delta_u$ with corresponding positions
$x''_{a,i}$, $y''_{a,i}$, $x''_{b,i}$ and $y''_{b,i}$ and
good event $\good''_1$ defined similarly as above,
that is,
\begin{eqnarray*}
&&\good''_1 = \{
\forall i \le L-1\dvtx
(i-1) \ell- 2M_t\delta_u < x''_{a,i},x''_{b,i} < (i-1) \ell,\\
&&\hspace*{141pt} i \ell<
y''_{a,i},y''_{b,i} < i \ell+ 2M_t\delta_u
\}.
\end{eqnarray*}

We define
\[
\good_1 = \good'_1\cup\good''_1.
\]
We show that this event
holds with high probability, conditioned on the sequence length
%
\begin{figure}

\includegraphics{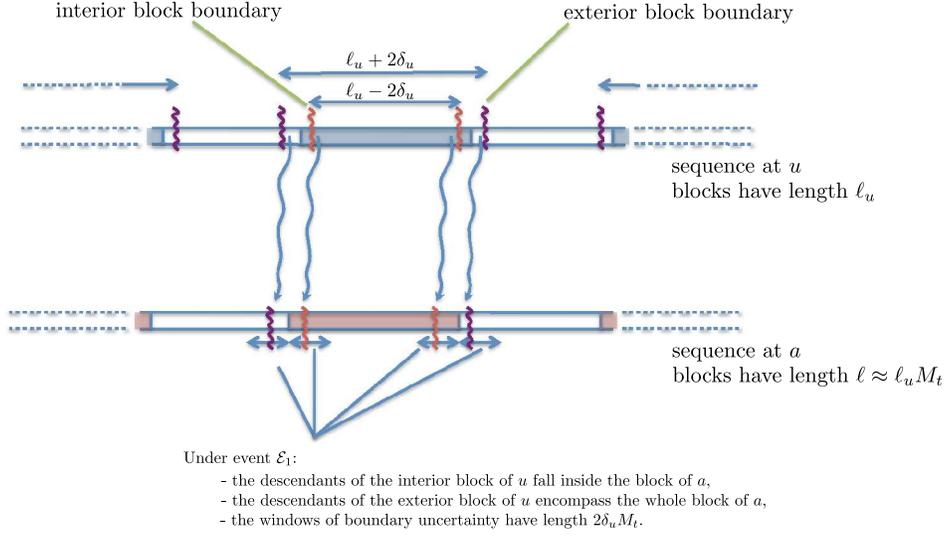}

\caption{Under the event $\mathcal{E}_1$ the descendants of the
interior blocks of $\sigma_u$ fall inside the corresponding blocks of
$\sigma_a$; the descendants of the exterior blocks of $\sigma_u$
contain all surviving sites inside the corresponding blocks of $\sigma
_a$; the windows of uncertainty have length $2M_t\delta_u$.}
\label{figgoodevent}
\end{figure}
$K_u$ at $u$ being at least $k_u'$. Figure~\ref{figgoodevent} shows
the structure of the indel process in the case that the event $\mathcal
{E}_1$ holds.
%
\begin{lemma}[(Interior/exterior block is inside/outside leaf
block)]\label{lemprobabilityofgoodevent1}
Conditioned on the event $\{K_u \ge k_u'\}$, we have
\[
\prob[\good_1] \ge1- {16 L \biggl(\frac{1}{k_u'}\biggr)^{\gamma}}.
\]
\end{lemma}
\begin{pf}
It follows from Lemma~\ref{lemmalengthconcentration}
that the left-most descendant of the site at position
$(i-1) \ell_u + \delta_u$ of $\s_u$ is located
inside the sequence of node $a$ at position at least
\begin{eqnarray*}
M_t \bigl((i-1) \ell_u + \delta_u\bigr) -
\scal_n\bigl( \sqrt{{k_u'} \log{k_u'}}\bigr)
&>& M_t \bigl((i-1) \ell_u + L \bigr)\\
&>& (i-1)\ell
\end{eqnarray*}
with probability $\ge1-(\frac{1}{k_u'})^{\gamma}$.
The other bounds follow similarly.
Taking a union bound over all $i$'s establishes the result.
\end{pf}

\subsubsection*{Block correlation}
Let $\js{ab,i}$ be the number of
common sites in the blocks $\s_{a,i}$ and
$\s_{b,i}$ that are
jointly surviving from $u$. Similarly, we define $\js{ab,i}'$
and $\js{ab,i}''$ where, for $\xi=a,b$, $\s_{\xi,i}'$ (resp.,
$\s_{\xi,i}''$) denotes the subsequence of $\s_\xi$ ranging from
position $x'_{\xi,i}$ (resp., $x''_{\xi,i}$)
to position $y'_{\xi,i}$ (resp., $y''_{\xi,i}$).
We define a good event for $\js{ab,i}$ as
\[
\good_2 = \{\forall i \leq L-1\dvtx\ell_u D_t^2 - 3 M_t \delta_u
\leq\js{ab,i}
\leq\ell_u D_t^2 + 3 M_t \delta_u\}.
\]

\begin{lemma}[(Jointly surviving sites in blocks)]
\label{lemprobabilityofgoodevent2}
Conditioned on the event $\{K_u \ge k_u'\}$, we have
\[
\prob[\good_2] \ge1- {18 L \biggl(\frac{1}{k_u'}\biggr)^{\gamma}}.
\]
\end{lemma}
\begin{pf}
We bound
\begin{eqnarray*}
\prob[\good_2^c] &=&
\prob[\good_2^c\cap\good_1] +
\prob[\good
_2^c\cap\good_1^c]
\leq
\prob[\good_2^c\cap\good_1] +
\prob[\good_1^c]\\
&\leq&
\prob[\good_2^c\cap\good_1] + 16 L \biggl(\frac
{1}{k_u'}\biggr)^{\gamma}.
\end{eqnarray*}
By construction, under $\good_1$ we have $\js{ab,i}' \leq\js{ab,i}
\leq\js{ab,i}''$
so that
\begin{eqnarray*}
\prob[\good_2^c\cap\good_1]
&\leq&
\prob[\exists i, \js{ab,i}' \leq\ell_u D_t^2 - 3 M_t \delta_u]\\
&&{} +
\prob[\exists i, \js{ab,i}'' \geq\ell_u D_t^2 + 3 M_t
\delta_u]\\
&\leq&
\prob\bigl[\exists i, \js{ab,i}' \leq(\ell_u - 2\delta_u +1
)D_t^2 -
\scal_n\bigl(\sqrt{k_u'\log k_u'}\bigr)
\bigr]\\
&&{} +
\prob\bigl[\exists i, \js{ab,i}'' \geq(\ell_u +
2\delta_u+1)D_t^2 +
\scal_n\bigl(\sqrt{k_u'\log k_u'}\bigr)
\bigr]\\
&\leq&
{2 L \biggl(\frac{1}{k_u'}\biggr)^{\gamma}}
\end{eqnarray*}
by Lemma~\ref{lemsurviving}, where we also used the fact that $D_t^2
\leq M_t$.
\end{pf}

\subsection{Estimation guarantees}

We are now ready to analyze the behavior of
our estimate $\epartdist(a,b)$. In this subsection
we compute the expectation and variance of
$\epartdist(a,b)$. We denote by $\realiz$
a realization of the indel process (but not of the substitution process)
on the paths between $u$ and $a,b$.
We denote by $\good$ the event that $\{K_u \geq k_u'\}$, $\good_1$ and
$\good_2$ are satisfied.
Suppose that $k_0 > k_0^*$ (defined in Section~\ref{secpartition}).
%
\begin{lemma}[(Block independence)] \label{lemindependence}
Conditioning on $\realiz$ and $\good$,
the \mbox{variables}
\[
\{\deviate_{a,2j+1}\deviate_{b,2j+1}\}_{j=1}^{L/2-1}
\]
are mutually independent.
\end{lemma}
\begin{pf}
Observe that when $K_u \ge k_u'$ the ancestral blocks $\s_{u,i}$
are well defined. Assuming that $k_0>k_0^*$, the interior blocks
$\s_{u,i}'$ are also well defined and disjoint. Hence, for a fixed
$\realiz$ under
$\good$, for all $i \le L-1$, both $\deviate_{a,i}$ and
$\deviate_{b,i}$ depend on the subsequence of $\s_u$ ranging
from position $(i-1) \ell_u - \delta_u +1$ to position $i \ell_u +
\delta_u-1$.
In this case, for $j \in\{1,\ldots,L/2-1\}$, different
$\deviate_{a,2j+1} \deviate_{b,2j+1}$'s are functions of different
subsequences
of $\s_u$. Observe that, since the root sequence is i.i.d. uniform and
the insertions above $u$ are also i.i.d. uniform,
the state of every site in $\s_u$ is uniform
and independent from the other sites. It follows from the above observations
that $\{\deviate_{a,2j+1} \deviate_{b,2j+1}\}_{j=1}^{L/2-1}$
are mutually independent.
\end{pf}
%
\begin{lemma}[(Expected correlation under good event)]\label{lemexpectation}
We have
\[
\expec[\deviate_{a,i} \deviate_{b,i} | \realiz,\good]
= \tfrac{1}{4} e^{-4\mutate t} e^{-2\death t} \ell_u
\pm
\scal_n\bigl(\sqrt{k_0 \log k_0}\bigr).
\]
\end{lemma}
\begin{pf}
Let $\deviate_{a,i}^\surv$ be the contribution to $\deviate_{a,i}$
from those common sites between $a$ and $b$ that are jointly
surviving from $u$. Let $\deviate_{a,i}^\nsurv= \deviate_{a,i} -
\deviate_{a,i}^\surv$,
and similarly for $b$. Then
\begin{eqnarray*}
\expec[\deviate_{a,i} \deviate_{b,i} | \realiz,\good]
&=& \expec[(\deviate_{a,i}^\surv+\deviate_{a,i}^\nsurv)
(\deviate_{b,i}^\surv+\deviate_{b,i}^\nsurv) | \realiz,\good]\\
&=& \expec[\deviate_{a,i}^\surv\deviate_{b,i}^\surv| \realiz
,\good],
\end{eqnarray*}
since the contribution from $\deviate_{a,i}^\nsurv$ and
$\deviate_{b,i}^\nsurv$ is independent and averages to $0$.
Write $\deviate_{a,i}^\surv$ as a sum over the jointly
surviving sites, that is,
\[
\deviate_{a,i}^\surv
= \sum_{j=1}^{\js{ab,i}} \biggl(z_{a,i}^{(j)} - \frac{1}{2}\biggr),
\]
where $z_{a,i}^{(j)}$ is $1$ if the corresponding site of $a$ is $0$.
Note that the terms in parentheses have zero expectation given $\realiz$
and $\good$.
Then,
\[
\expec[\deviate_{a,i}^\surv\deviate_{b,i}^\surv| \realiz,\good]
= \sum_{j=1}^{\js{ab,i}} \expec\biggl[\biggl(z_{a,i}^{(j)} - \frac
{1}{2}\biggr)\biggl(z_{b,i}^{(j)} - \frac{1}{2}\biggr) \Big| \realiz
,\good\biggr]
\]
by independence of the sites.
We compute the expectation above.
We have
\begin{eqnarray*}
\expec\biggl[\biggl(z_{a,i}^{(j)} - \frac{1}{2}\biggr)
\biggl(z_{b,i}^{(j)} - \frac{1}{2}\biggr)\Big | \realiz,\good\biggr]
&= &\expec\biggl[\biggl(z_{a,i}^{(j)}z_{b,i}^{(j)} - \frac
{1}{2}z_{a,i}^{(j)} - \frac{1}{2}z_{b,i}^{(j)} + \frac{1}{4}
\biggr) \Big| \realiz,\good\biggr]\\
& = &\expec\bigl[z_{a,i}^{(j)}z_{b,i}^{(j)} | \realiz,\good
\bigr] - \frac{1}{4}\\
& = &\frac{1}{2}\frac{1+e^{-4\mutate t}}{2} - \frac{1}{4}\\
& = &\frac{1}{4}e^{-4\mutate t}.
\end{eqnarray*}
Therefore,
\[
\expec[\deviate_{a,i}^\surv\deviate_{b,i}^\surv| \realiz,\good]
= \tfrac{1}{4}e^{-4\mutate t} \js{ab,i}.
\]
The result then follows from the definition of $\good_2$.
\end{pf}
%
\begin{lemma}[(Variance under good event)]\label{lemvariance}
We have
\[
\var[\deviate_{a,i} \deviate_{b,i} | \realiz,\good] \le\tfrac
{3}{16} \ell^2.
\]
\end{lemma}
\begin{pf}
By Cauchy--Schwarz we have
\begin{eqnarray*}
\expec[\deviate_{a,i}^2 \deviate_{b,i}^2 | \realiz,\good]
&\leq& (\expec[\deviate_{a,i}^4 | \realiz,\good] \expec[\deviate
_{b,i}^4 | \realiz,\good])^{1/2}\\
&=& \bigl(\tfrac{1}{16}(3\ell^2 - 2\ell) \cdot
\tfrac{1}{16}(3\ell^2 - 2\ell)\bigr)^{1/2}\\
&\leq& \tfrac{3}{16} \ell^2,
\end{eqnarray*}
where we used the fact that the length of the sequences $\s_{a,i}$ and
$\s_{b,i}$ is deterministically $\ell$, and the number of zeros in
$\s_{a,i}$ and $\s_{b,i}$ follows a binomial distribution with $\ell
$ trials and probability $1/2$.
\end{pf}
%
\begin{lemma}[(Distance estimate)]\label{lemexpectation-variance}
We have
\[
\expec[\epartdist(a,b) | \realiz,\good]
= \tfrac{1}{4} e^{-(4\mutate+ \death+ \birth) t} \ell
\pm
\scal_n\bigl(\sqrt{k_0 \log k_0} \bigr)
\]
and
\[
\var[\epartdist(a,b) | \realiz,\good]
\le\frac{3}{8} {1 \over\lfloor k_0^{1-\zeta} \rfloor} \ell^2.
\]
In particular, the standard deviation
\[
\operatorname{STD}[\epartdist(a,b) | \realiz,\good]
= O\bigl(k_0^{({3\zeta- 1})/{2}}\bigr) = o\bigl(\sqrt{k_0}\bigr)
\]
for $\zeta> 1/2$ small enough.
\end{lemma}
\begin{pf}
From Lemma~\ref{lemindependence}, the $L/2 = \lfloor k_0 / \ell
\rfloor/2$ terms in $\epartdist(a,b)$
are mutually independent.
The proof then follows from Lemmas~\ref{lemexpectation} and \ref
{lemvariance}
and the definition of $\ell_u$.
\end{pf}

\subsection{Concentration}

We now show that our distance estimate is concentrated. For notational
convenience, we denote by $\prob_u^*$ the probability measure induced
by conditioning on the event $\{K_u \geq k_u'\}$. Recall that the event
$\mathcal{E}$ is contained in \mbox{$\{K_u \geq k_u'\}$}.
%
\begin{lemma}[(Concentration of distance estimate)]
\label{lemdistanceconcentration1}
Let $\alpha> 0$ be such that
$\zeta- \alpha> 1/2$,
and
$\beta= 1 - \zeta- 2\alpha> 0$
for $\zeta> 1/2$ small enough.
Then for $k_0$ large enough
\[
\prob_u^*\biggl[\biggl|\frac{4}{\ell}\epartdist(a,b)
- e^{-(4\mutate+ \death+ \birth) t} \biggr|
> \frac{1}{k_0^{\alpha}} \biggr]
\leq O\biggl(\frac{1}{k_0^{\beta}}\biggr).
\]
\end{lemma}
\begin{pf}
We use Chebyshev's inequality. We first condition on $\realiz,\good$.
Recalling that $\ell= \lceil k_0^{\zeta} \rceil$, note that
\begin{eqnarray*}
&&\prob_u^*\biggl[\frac{4}{\ell}\epartdist(a,b)
> e^{-(4\mutate+ \death+ \birth) t} + \frac{1}{k_0^{\alpha
}} \Big| \realiz,\good\biggr]\\[-2pt]
&&\qquad \leq\prob_u^*\biggl[\epartdist(a,b)
> \frac{\ell}{4}e^{-(4\mutate+ \death+ \birth) t}
+ \frac{\ell}{4}\frac{1}{k_0^{\alpha}} \Big| \realiz,\good\biggr]\\[-2pt]
&&\qquad \leq\prob_u^*\biggl[\epartdist(a,b)
> \expec[\epartdist(a,b) | \realiz,\good] -
\scal_n\bigl(\sqrt{k_0 \log k_0} \bigr)
+ \frac{\ell}{4}\frac{1}{k_0^{\alpha}} \Big| \realiz,\good\biggr]\\[-2pt]
&&\qquad \leq\frac{{3}\ell^2 /({8}{\lfloor k_0^{1-\zeta}
\rfloor}) }
{({\ell}/({4}{k_0^{\alpha}}) -
\scal_n(\sqrt{k_0 \log k_0} )
)^2}\\[-2pt]
&&\qquad = O\biggl(\frac{1}{k_0^{1-\zeta-2\alpha}}\biggr).
\end{eqnarray*}
The other direction is similar. Taking expectation over $\realiz$, we have
\[
\prob_u^*\biggl[\biggl|\frac{4}{\ell}\epartdist(a,b)
- e^{-(4\mutate+ \death+ \birth) t} \biggr|
> \frac{1}{k_0^{\alpha}} \Big| \good\biggr]
\leq O\biggl(\frac{1}{k_0^{\beta}}\biggr).
\]
Choose $\gamma>0$
in Lemmas~\ref{lemmalengthconcentration} and~\ref{lemsurviving}
large enough so that
\[
\gamma- (1 -\zeta) > \beta.
\]
Then, from Lemmas~\ref{lemprobabilityofgoodevent1}
and~\ref{lemprobabilityofgoodevent2},
we have
\begin{eqnarray*}
&&\prob_u^*\biggl[\biggl|\frac{4}{\ell}\epartdist(a,b)
- e^{-(4\mutate+ \death+ \birth) t} \biggr|
> \frac{1}{k_0^{\alpha}}\biggr]\\[-2pt]
&&\qquad \leq
\prob_u^*\biggl[\biggl|\frac{4}{\ell}\epartdist(a,b)
- e^{-(4\mutate+ \death+ \birth) t} \biggr|
> \frac{1}{k_0^{\alpha}} \Big| \good\biggr]\prob_u^*[\good]
+\prob_u^*[\good^c]\\[-2pt]
&&\qquad \leq O\biggl(\frac{1}{k_0^{\beta}}\biggr).
\end{eqnarray*}
\upqed\end{pf}

The proofs of Theorems~\ref{thmconsistency} and~\ref{thmmain} are
given in the next section.

\section{Putting it all together}
\label{secalgo}

\subsection*{Large-scale asymptotics}
%
We are ready to prove our main result in the molecular clock case. We
postpone the more general case to the next section.
A last bit of notation:
for a pair of leaves $a, b \in[n]$, we denote by $t_{ab}$ the time
between $a$, $b$ and their most recent common ancestor.
%
\begin{pf*}{Proof of Theorem~\ref{thmmain}}
We first give a bound on the diameter of the tree.
Let $h$ (resp., $H$) be the length of the shortest
(resp., longest) path between the root and a leaf in graph
distance. Because the number of leaves is $n$ we must have
$2^h \leq n$ and $2^H \geq n$. Since all leaves are contemporaneous it
must be that
$Hf \leq hg$. Combining these constraints gives that the diameter
$\diam$ satisfies
\[
2 \frac{f}{g} \log_2 n \leq2 h \leq\diam\leq2H \leq2\frac{g}{f}
\log_2 n.
\]
Given our bound on the diameter of the tree, it follows that the time
from the root $r$ of the tree to any leaf is at most $\frac{g^2}{f}
\log_2 n$. Suppose that the length $k_r$ at the root of the tree
satisfies $k_r > k_r^* = k_r^*(k_0)$, where $k_r^*$ is the minimum
integer satisfying
\[
k_r^* \ge e^{({ g^2}/{f}) \cdot{\death} \log_2 n} \bigl(k_0 +
\scal_n\bigl(\sqrt{k_r^* \log{k_r^*}}\bigr)
\bigr),
\]
where the small polynomial factor is taken to be the one
used in Lemma~\ref{lemmalengthconcentration}.

Lemma~\ref{lemmalengthconcentration} and a union bound then imply
that with probability at least
\[
1- O(n) \cdot( {k^*_r} )^{-\gamma}
\]
for all nodes $u$
\[
K_u \ge k_u'.
\]

\begin{lemma}[(Concentration of distance estimate)]\label{lemalldistances}
For all $\alpha' > 0$, $\beta'>0$, there exists
$k_0 = n^{\beta'''}$ with $\beta''' > 0$
large enough so that if the sequence length at the root is $k_r >
k_r^*(k_0)$, then
\[
\prob\biggl[\forall a,b\in[n], \biggl|\frac{4}{\ell}\epartdist(a,b)
- e^{-(4\mutate+ \death+ \birth) t_{ab}} \biggr|
\leq\frac{1}{n^{\alpha'}} \biggr]
=1-O\biggl({1 \over n^{\beta'}}\biggr).
\]
\end{lemma}
\begin{pf}
This follows from Lemma~\ref{lemdistanceconcentration1} and our
observation above that, if $k_r > k_r^*(k_0)$, with probability at
least $1- O(n) \cdot( {k^*_r} )^{-\gamma}$, then $K_u
\ge k_u'$ for all nodes~$u$.
\end{pf}

Given our bound on the diameter of the tree, it follows
that for all pairs of leaves $a$, $b$ and small $\eps> 0$
\[
e^{-(4\mutate+ \death+ \birth) t_{ab} \pm\eps} = e^{-(4\mutate+
\death+ \birth) t_{ab}} \bigl(1\pm O(\eps)\bigr) \ge\frac{1}{n^{\alpha
''}}\bigl(1\pm O(\eps)\bigr).
\]
Therefore, choosing $\alpha'$ large enough in Lemma~\ref{lemalldistances},
we get that all distances can be estimated within a small $\eps$ simultaneously
with probability going to~$1$.

Using the standard Buneman algorithm, we can recover the tree efficiently
(see, e.g.,~\cite{SempleSteel03}).
\end{pf*}

\subsection*{Constant-size case}
The proof of Theorem~\ref{thmconsistency} for the molecular clock
case builds on the proof of Theorem~\ref{thmmain} by treating $n$ as
a constant and letting the sequence length at the root of the tree go
to infinity.
%
\begin{pf*}{Proof of Theorem~\ref{thmconsistency} \textup{(Molecular
clock case)}}
We can restate Lem\-ma~\ref{lemalldistances} in the following form,
where the failure probability is expressed more cleanly in terms of the
sequence length at the root of the tree. The proof of the lemma is
essentially the same.
%
\begin{lemma}[(Concentration of distance estimate)]
\label{lemalldistancesrestated}
For all $\alpha' > 0$, there exists
$k_0^* = n^{\beta'''}$ for $\beta''' > 0$
large enough such that if the sequence length at the root is $k_r >
k_r^*(k_0^*)$, then
\begin{eqnarray*}
&&\prob\biggl[\forall a,b\in[n], \biggl|\frac{4}{\ell}\epartdist(a,b)
- e^{-(4\mutate+ \death+ \birth) t_{ab}} \biggr|
\leq\frac{1}{n^{\alpha'}} \biggr]\\
&&\qquad=1- O(n \cdot k_r^{-\gamma})-O(n^2 \cdot k_r^{-\beta}).
\end{eqnarray*}
\end{lemma}

Repeating the proof of Theorem~\ref{thmmain} above, it follows that
the algorithm fails to reconstruct the phylogeny with probability $O(n
\cdot k_r^{-\gamma})+O(n^2 \cdot k_r^{-\beta})$. Letting $k_r
\rightarrow+\infty$ concludes the proof of Theorem~\ref{thmconsistency}.
\end{pf*}

\section{Extensions} \label{secextensions}

\subsection*{GTR model}
We briefly discuss how our results can be extended
to GTR models. For background on GTR models, see, for example, \cite
{Felsenstein04}.
Let $Q$ be a reversible $4\times4$ rate matrix with stationary
distribution $\pi$. Our new sequence evolution process is identical
to the one described in Definition~\ref{defevolutionaryprocessbranch}
except that the substitution process
is a continuous-time Markov process with rate matrix $\mutate_e Q$.
The rate matrix $Q$ has 4 nonnegative eigenvalues.
For convenience, we assume that the largest negative eigenvalue
is $-1$. We denote by $w$ the corresponding eigenvector which
we assume is normalized as
\[
\sum_{s \in\{\A,\G,\C,\T\}} \pi_s w_s^2 = 1.
\]
We now perform the following transformation of the state space.
For a node~$u$, let $\s_u = (\s_u^1,\ldots,\s_u^{K_u})$ be the transformed
sequence at $u$ where $\s_u^i = w_\A$ (resp., $w_\G,w_\C,w_\T$) if
the state at site $i$ is $\A$ (resp., $\G,\C,\T$). Note that, under
stationarity,
the expectation of the state at site $i$ is $0$ by orthogonality of
$\pi$ and $w$.
Then our distance estimator is
\[
\edist(a,b) = \Biggl(\sum_{i=1}^{K_a}\s_a^i\Biggr)
\Biggl(\sum_{j=1}^{K_b}\s_b^j\Biggr).
\]
In particular, in the two-state CFN case, we have $w = (+1,-1)$
and we obtain the same estimate as before, up to a constant.
We now analyze the expectation of this quantity.
For $u \in V$, we let
\[
\deviate_u = \sum_{i=1}^{K_u}\s_u^i.
\]

\begin{lemma}\label{lem1prime}
The following holds:
%
\begin{equation}\label{eqsingleprime}
\expec[\deviate_a | \s_r] = e^{-(\mutate+\death) t} \deviate_r.
\end{equation}
\end{lemma}
%
\begin{remark}\label{remnormalize}
Note that this formula is slightly different than that in
Lem\-ma~\ref{lem1} because of the normalization implied by
requiring $Q$ to have second eigenvalue $-1$.
\end{remark}
\begin{pf*}{Proof of Lemma~\ref{lem1prime}}
The sites created after $r$ contribute $0$ in expectation.
Of course, so do the deleted sites. The fraction of sites that
survive is $e^{-\death t}$. Suppose site $i$ survives, then note that
\[
\expec[\s_a^i | \s_r^i = w_s, \mbox{$i$ survives}]
= \sum_{s'}(e^{\mutate t Q})_{s s'} w_{s'} = e^{-\mutate t} w_s.
\]
Summing over all sites of $r$ we get
\[
\expec[\deviate_a | \s_r] = e^{-(\mutate+\death) t} \deviate_r
\]
as claimed.
\end{pf*}

Consider now a ``fork'' tree, that is,
a root $r$ from which emanates a single edge $e_u = (r,u)$ which in turn
branches into two edges $e_a = (u,a)$ and $e_b = (u,b)$.
For $x=a,b,u$, we denote the parameters of edge $e_x$ by
$t_x, \birth_x, \death_x, \mutate_x$.
Our goal is to compute $\expec[\edist(a,b)]$ assuming that the
sequence length
at the root is $k$. The proof is similar to Lemma~\ref{lem2}.
%
\begin{lemma}\label{lem2prime}
The following holds:
\[
\expec[\edist(a,b)]
= e^{-(\mutate_a+\death_a) t_a}
e^{-(\mutate_b+\death_b) t_b}
e^{-(\death_u - \birth_u)t_u} k.
\]
Note that Remark~\ref{remnormalize} also applies here.
\end{lemma}
\begin{pf}
We have
\begin{eqnarray*}
\expec[\edist(a,b)]
&=& \expec[\deviate_a \deviate_b]\\
&=& \expec[\expec[\deviate_a \deviate_b | \s_u
]]\\
&=& \expec[\expec[\deviate_a | \s_u]\expec
[\deviate_b | \s_u]]\\
&=& e^{-\mutate_a t_a} e^{-\death_a t_a} e^{-\mutate_b t_b}
e^{-\death_b t_b} \expec[\deviate_u^2]\\
&=& e^{-\mutate_a t_a} e^{-\death_a t_a} e^{-\mutate_b t_b}
e^{-\death_b t_b} \expec[\expec[\deviate_u^2 |
K_u]]\\
&=& e^{-\mutate_a t_a} e^{-\death_a t_a} e^{-\mutate_b t_b}
e^{-\death_b t_b} \expec[\var[\deviate_u | K_u
]]\\
&=& e^{-\mutate_a t_a} e^{-\death_a t_a} e^{-\mutate_b t_b}
e^{-\death_b t_b} \expec[K_u \var[\s_u^1]]\\
&=& e^{-\mutate_a t_a} e^{-\death_a t_a} e^{-\mutate_b t_b}
e^{-\death_b t_b} \expec[K_u \expec[(\s_u^1)^2]]\\
&=& e^{-\mutate_a t_a} e^{-\death_a t_a}
e^{-\mutate_b t_b} e^{-\death_b t_b}
e^{-(\death_u - \birth_u)t_u} k
\end{eqnarray*}
by Lemma~\ref{lem1prime}.
\end{pf}

From the previous lemmas, one can adapt the proofs above
to the GTR case.

\subsection*{Nonclock case}
Using Lemma~\ref{lem2prime}, we can get rid of the molecular clock
assumption. Consider again the fork tree, but assume that each
edge is in fact a path. An adaptation of Lemma~\ref{lem2prime}
gives the following lemma.
%
\begin{lemma}\label{lem3prime}
The following holds:
\[
-\ln\biggl(\frac{\expec[\edist(a,b)]}{\sqrt{\expec
[K_a] \expec[K_b]}}\biggr)
= \sum_{e\in\path(a,b)}(\mutate_e+\death_e/2+\birth_e/2) t_e.
\]
Note that Remark~\ref{remnormalize} also applies here.
\end{lemma}
\begin{pf}
Note that
\[
-\ln(k^{-1}\expec[K_a]) = \sum_{e\in\path(r,a)}
(\death_e - \birth_e)t_e
\]
and similarly for $b$. A variant of Lemma~\ref{lem2prime} gives
\[
-\ln(k^{-1}\expec[\edist(a,b)])
= \sum_{e\in\path(a,b)}(\mutate_e+\death_e) t_e
+ \sum_{e\in\path(r,u)} (\death_e - \birth_e)t_e.
\]
The result follows by subtracting the previous expressions.
\end{pf}

The expression in Lemma~\ref{lem3prime} provides the additive
metric needed to extend our results to nonclock
bounded-rates case.

\section{Concluding remarks}

We have shown how to reconstruct phylogenies
under the bounded-rates, GTR model with
indels.
Our efficient algorithm requires polynomial-length
sequences at the root. A natural open
problem arises from this work:
Can our results be extended to general
trees with bounded branch lengths, as opposed
to the bounded-rates model? The key difference
between the two models is that the former
may have a linear diameter whereas the latter
has logarithmic diameter. To extend our results,
one would need to deal with far away
leaves that are almost uncorrelated but for which
our block structure does not apply.
%



\printaddresses

\end{document}